\documentclass{article}  
\usepackage{amssymb}              
\usepackage{amsthm}
\usepackage{amsmath}
\usepackage{eucal}
\usepackage{verbatim}
\usepackage{graphicx}
\usepackage{multirow}
\usepackage{fancyhdr}
\usepackage{color}
\usepackage{enumerate}
\usepackage{calrsfs}
\usepackage{cancel}
\usepackage{fullpage}
\usepackage{amsmath}

\newtheorem{theorem}{Theorem}

\newtheorem{definition}{Definition}

\newtheorem{proposition}{Proposition}

\def\R{\mathbb{R}}

\newcommand\Vset{V}
\newcommand\Eset{E}
\newcommand\En{\mathcal{E}}

\usepackage{color}
\definecolor{lightgray}{gray}{0.75}

\newcommand\greybox[1]{%
  \vskip\baselineskip%
  \par\noindent\colorbox{lightgray}{%
    \begin{minipage}{\textwidth}#1\end{minipage}%
  }%
  \vskip\baselineskip%
}

\begin{document}

\title{Numerical approach to centrality of optimal transportation networks}

\date{}

\pagenumbering{arabic}

\author{Jan Haskovec\footnote{Mathematical and Computer Sciences
            and Engineering Division,
         King Abdullah University of Science and Technology,
         Thuwal 23955-6900, Kingdom of Saudi Arabia;
         {\it jan.haskovec@kaust.edu.sa}}}

\maketitle
\vspace{2mm}

\textbf{Abstract.}
We study hierarchical properties of optimal transportation networks with biological background.
The networks are obtained as minimizers of an energy functional which involves a metabolic cost term of a power-law form with exponent $\gamma>0$.
In the range $\gamma\in (0,1)$, most relevant for biological applications, the functional is non-convex and its local minima correspond to loop-free graphs (trees).
We propose a numerical scheme that performs energy descent by searching the discrete set of local minimizers, combined with a Monte-Carlo approach.
We verify the performance of the scheme in the borderline case $\gamma=1$, where the functional is convex.
For~a~particular example of a leaf-shaped planar graph, we evaluate the global reaching centrality (GRC) of the (local) minimizers in dependence on the value of $\gamma\in (0,1]$.
We observe that the GRC, which can be understood as a measure of hierarchical organization of the graph, monotonically increases with increasing $\gamma$.
To our best knowledge, this is the first quantification of the influence of the value of the metabolic exponent on the hierarchical organization of the (almost) optimal transportation network.
\vspace{2mm}

\section{The discrete network formation model}\label{sec:1}
In this paper we focus on the discrete network formation model introduced by Hu and Cai \cite{Hu-Cai-13, Hu-Cai-19} and further studied, including its various modifications
and extensions, analytically in, e.g., \cite{BHMR, HKM, HMP15, HMP22, HMP23} and numerically in, e.g., \cite{bookchapter, Ast1, Ast2, Di, HMPS}.
Let $G=(\Vset,\Eset)$ denote a prescribed undirected connected graph with a finite set of vertices $\Vset$ and a finite set of edges $\Eset\subseteq\Vset\times\Vset$,
where each edge $(i,j) \equiv (j,i) \in \Eset$ represents an undirected connection between vertices $i\neq j\in\Vset$.
We assume that the graph $G$ is connected, i.e., that every pair of vertices in $\Vset$ is connected by at least one path (sequence of edges) in $\Eset$.
The number of vertices shall be denoted by  $|\Vset|$ in the sequel.
For each edge $(i,j)\in\Eset$ of the graph $G$ we prescribe its length $L_{ij}=L_{ji}>0$, and, moreover,
at each vertex $i \in\Vset$ we prescribe the strength of the external material flow $S_i\in\R$, where we adopt the convention that $S_i$ is positive for inflow nodes (sources) and negative for sinks. 

The goal is to optimize the transportation network in terms of the edge conductivities $C_{ij} \geq 0$, for all $(i, j)\in\Eset$.
In order to define the energy functional, we introduce the material flow pressures $P_i \in \R$ for $i\in\Vset$.
Assuming that the flow takes place in the laminar (Poiseuille) regime, the oriented material flux from vertex $i$ to vertex $j$ is given by
\begin{align}\label{eq:flowrate}
   Q_{ij} := C_{ij}\frac{P_i-P_j}{L_{ij}} \qquad\text{for all~}(i,j)\in \Eset.
\end{align}
Moreover, we impose the local mass conservation in each vertex, expressed in terms of the Kirchhoff law
\begin{align}\label{eq:kirchhoff}
   \sum_{j\in N(i)} C_{ij}\frac{ P_i-P_j}{L_{ij}}=S_i\qquad \text{for all~}i\in \Vset.
\end{align}
Here $N(i)$ denotes the set of vertices connected to $i\in\Vset$ through an edge, i.e.,
$N(i) := \{ j\in\Vset;\, (i,j)\in\Eset \}.$
Clearly, a necessary condition for the solvability of \eqref{eq:kirchhoff} is the global mass conservation
 $\sum_{i\in\Vset} S_i = 0$, which we assume in the sequel.
 For a detailed discussion of the solvability of the linear system of equations \eqref{eq:kirchhoff} for the vector of pressures $P=(P_i)_{i\in\Vset}$
we refer to \cite[Remark 1]{HV23}.

The energy expenditure of the transportation network is, according to the model \cite{Hu-Cai-13, Hu-Cai-19}, given by
\begin{align}\label{eq:energy}
   \En[C] := \sum_{(i,j)\in\Eset}\left( \frac{Q_{ij}[C]^2}{C_{ij}}+\frac{\nu}{\gamma} C_{ij}^{\gamma}\right) L_{ij},
\end{align}
where $Q_{ij}=Q_{ij}[C]$ is given by \eqref{eq:flowrate} with pressures calculated from \eqref{eq:kirchhoff},
and $\nu>0$ is the so-called metabolic coefficient.
The value of the exponent $\gamma>0$ depends on the specific biological application.
For instance, $\gamma\in [1/2, 1]$ for models of leaf venation \cite{Hu-Cai-13, Laguna-Bohn-Jagla, Runions} and 
$\gamma = \frac{1}{2}$ in blood vessel systems \cite{Murray}.
In this paper we shall consider \eqref{eq:energy} with $\gamma\in (0,1]$.

Let us note that the first term in \eqref{eq:energy} is convex in $C\in\R_+^{|V|}$, see \cite[Lemma 1]{HV23}, while
for $\gamma\in(0,1)$ the second term is concave. Consequently, for $\gamma\in(0,1)$
the functional $\En=\En[C]$ is non-convex.
Moreover, combining the results of \cite[Theorem 2.1]{BHMR} and \cite[Theorem 1]{HV23},
we have the following characterization of the local minimizers of the constrained minimization problem \eqref{eq:kirchhoff}--\eqref{eq:energy}.

\begin{proposition}\label{prop:1}
Let $\gamma\in (0,1)$. Then every local minimizer $C\in\R^{|V|}_+$ of the constrained problem \eqref{eq:kirchhoff}--\eqref{eq:energy}
is loop-free, where we treat edges $(i,j)\in\Eset$ with $C_{ij}=0$ as nonexistent.
Conversely, every spanning tree of the connected graph $G=(V,E)$ induces a local minimizer,
i.e., there exists a set of conductivities $C\in\R^{|V|}_+$, supported on the spanning tree, which is a local minimizer
of \eqref{eq:kirchhoff}--\eqref{eq:energy}.
\end{proposition}

It follows that in order to find a global minimizer of \eqref{eq:kirchhoff}--\eqref{eq:energy}, it suffices to search
the set of spanning trees of the graph $G=(V,E)$. However, this task is computationally infeasible even for graphs
of moderate size. For instance, the number of spanning trees for a complete graph with $|\Vset|$ nodes is $|\Vset|^{|\Vset|-2}$;
for planar graphs, the number of possible spanning trees grows exponentially, see, e.g., \cite{Buchin-Schulz}.
Consequently, in this paper we propose and implement a simple energy descent algorithm which
in each step randomly selects one edge from the spanning tree and seeks a replacement by another edge from $\Eset$
such that the resulting graph has a lower value of the energy. This procedure, described in detail in Section \ref{sec:alg}
below, of course does not guarantee the global minimizer to be found. However, as we shall demonstrate on an example
in Section \ref{sec:results}, it likely provides solutions that are reasonably close to the optimum.
In particular, one can test the method rigorously against known global minimizers of the energy $\En=\En[C]$ with $\gamma=1$.
This is due to the fact that for $\gamma=1$ the energy is convex and differentiable in $C\in\R_+^{|V|}$,
so that global minimizers can be reliably found, e.g., by gradient descent methods.
Moreover, \cite[Theorem 2]{HV23} provides the following characterization of the set of minimizers of $\En=\En[C]$
for $\gamma=1$.

\begin{proposition}\label{prop:2}
Let $\gamma = 1$. Then the set $M$ of minimizers of \eqref{eq:kirchhoff}--\eqref{eq:energy}
is a nonempty closed convex subset of $\R_+^{|V|}$.
Its extremal points represent loop-free graphs, and vice-versa, the loop-free elements of $M$ are extremal points of $M$.
\end{proposition}

Therefore, our method of searching the set of spanning trees can be applied also for the case $\gamma=1$,
and the (local) minimum found can be compared with the global minimum provided by
the gradient descent algorithm. We shall carry out this test in Section \ref{sec:results}
for a particular example of a planar graph and observe that our heuristical method
is indeed able to discover the global minimizer in this case.

Finally, having reliable numerical results at our disposal, we investigate how the hierarchical structure
of the (almost) optimal transportation networks depends on the value of the metabolic exponent $\gamma\in (0,1]$.
We use the global reaching centrality (GRC) proposed in \cite{Vicsek} as a measure of the degree of hierarchy.
As we shall see in Section \ref{sec:results}, the value of the GRC coefficient increases with increasing value of $\gamma\in (0,1]$.

\section{Discrete energy descent on trees for $\gamma<1$}\label{sec:alg}

If the set of "active" edges, i.e., edges $(i,j)\in\Eset$ with $C_{ij}>0$, forms
a spanning tree of $G=(\Vset, \Eset)$, then the fluxes $Q_{ij}$
through the active edges are determined uniquely by the sources/sinks $S_i$.
Indeed, for any active edge $(i,j)\in\Eset$, there is a unique division of the vertices into two disjoint sets $V^{(i)}$, $V^{(j)}$
such that nodes from $V^{(i)}$ are connected to $i$ by paths in $E \setminus \{(i,j)\}$, and analogously for nodes in $V^{(j)}$.
By the mass conservation, the flux through $(i,j)$ is then given by
\begin{equation}  \label{eq:Qij}
   Q_{ij} = \sum_{k\in V^{(i)}} S_k - \sum_{k\in V^{(j)}} S_k.
\end{equation}
With the explicit formula \cite[Lemma 2.1]{HKM} for the derivative of \eqref{eq:energy},
\[
   \frac{\partial \En[C]}{\partial C_{ij}} = \left( - \frac{Q_{ij}^2}{C^2_{ij}}+ \nu C_{ij}^{\gamma-1}\right) L_{ij},
\]
the optimal conductivities are given by $C_{ij}:= \left( Q_{ij}^2/\nu \right)^{1/(\gamma +1)}$,
which covers also the inactive edges with $Q_{ij} = C_{ij} = 0$.
The value for of the optimal energy $\En^\ast = \En^\ast[Q]$ is then
\begin{equation}  \label{eq:Enast}
   \En^\ast[Q] =  2 \nu^\frac{\gamma}{\gamma+1}  \sum_{(i,j)\in\Eset}  |Q_{ij}|^\frac{2\gamma}{\gamma+1} L_{ij}.
\end{equation}

Our discrete descent algorithm for minimization of \eqref{eq:kirchhoff}--\eqref{eq:energy} with $\gamma\in (0,1)$
is based on searching the set of spanning trees of the graph $G=(\Vset,\Eset)$,
relying on Proposition \ref{prop:1}.

\greybox{Discrete descent algorithm on trees\hfill}
\vspace{-5mm}

\begin{itemize}
\item
Initialize $T$ as a spanning tree (randomly) extracted 
from $G=(\Vset, \Eset)$
\item
Repeat:
\begin{itemize}
\item
Select an edge $(I,J)$ from $T$ randomly
\item
For all edges $(i,j)\in\Eset$:
\begin{itemize}
\item
Calculate fluxes $Q^{(i,j)}$, using formula \eqref{eq:Qij}, for the tree $\tilde T^{(i,j)}$ obtained from
 $T$ by removing edge $(I,J)$ and adding edge $(i,j)$
\item
Calculate the energy $\En^\ast[Q^{(i,j)}]$ for $\tilde T^{(i,j)}$, given by formula \eqref{eq:Enast}
\item
Record the edge $(i,j)$ that gives the lowest value of $\En^\ast[Q^{(i,j)}]$
\end{itemize}
\item
If $\min_{(i,j)\in\Eset} \En^\ast[Q^{(i,j)}]$ is lower than the energy of $T$, update $T$ by deleting edge $(I,J)$ and adding edge $(i,j)$
\end{itemize}
\item
If all edges $(I,J)$ in $T$ have been tried and no energy decreasing replacement $(I,J) \to (i,j)$ has been found,
stop and output $T$
\end{itemize}
\vspace{-5mm}

\noindent\hrulefill
\vspace{3mm}

In summary, in each iteration the discrete descent algorithm tries to find
a tree with lower value of the energy by replacing one edge of the tree $T$
by another edge from the underlying graph $G$. This is repeated as long as
such a pair of edges is found. Obviously, if the descent of the energy is not
achievable any more by swapping a single pair of edges, one could start searching
for $k$ pairs of edges to be swapped, taking progressively $k=2,3,\dots$. This would obviously
lead to exponentially increasing computational cost in $k$.
Instead, we chose to repeat the runs of the algorithm with $k=1$ in a Monte-Carlo approach,
each time initializing with a randomly generated spanning tree $T$ from $G=(\Vset, \Eset)$.
Among all outputs of the independent Monte-Carlo runs,
we select the tree with the lowest value of the energy as the final output.
This way we were able to achieve results that seem to be
sufficiently close to the global minimum, while keeping the computational cost low.
Moreover, as the Monte Carlo runs are independent from each other,
this approach is perfectly parallelizable.

In fact, one can rigorously test the performance
of the discrete descent algorithm in the borderline case $\gamma=1$.
Then, according to Proposition \ref{prop:2}, the global minimizer of the energy is a convex set
containing a tree.
A global minimizer can be reliably found by a gradient descent scheme
since with $\gamma=1$ the energy functional is convex and differentiable \cite[Lemma 2.1]{HKM}.
We carried out this test for the leaf-shaped graph described in Section \ref{sec:results} below,
performing $1000$ runs of the Monte-Carlo simulation.
The discrete descent algorithm found the exact global minimizer in $40$ runs,
i.e., $4\%$ of the time.
Moreover, in $99\%$ percent of the runs, the energy of the graph
found by the algorithm lied within $1\%$ relative error
with respect to the global energy minimum.

\section{Global reaching centrality for a leaf-shaped example}\label{sec:results}

\begin{figure}[h] \centering
 \includegraphics[width=.5\textwidth]{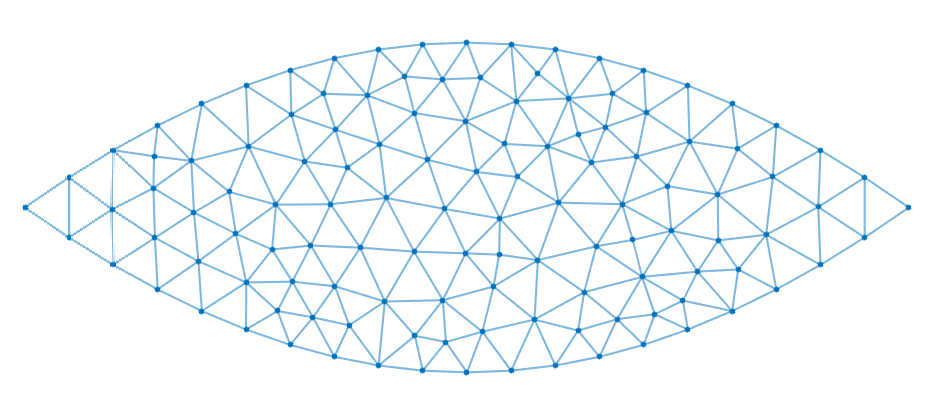}
 \caption{The planar leaf-shaped graph with 122 nodes and 323 edges, introduced in \cite[Fig. 5]{HV23} \label{fig:leaf}.}
\end{figure}

We provide outputs of the discrete descent algorithm presented in Section \ref{sec:alg},
applied to a leaf-shaped graph introduced in \cite{HV23}.
The planar graph consists of $|V|=122$ nodes and $323$ edges, see Fig. \ref{fig:leaf}.
There is a single source $S_i=1$ in the left-most node
(``stem" of the leaf), while $S_j = -(|V|-1)^{-1}$ for all other nodes.

\begin{figure}[h!] \centering
 \includegraphics[width=.45\textwidth]{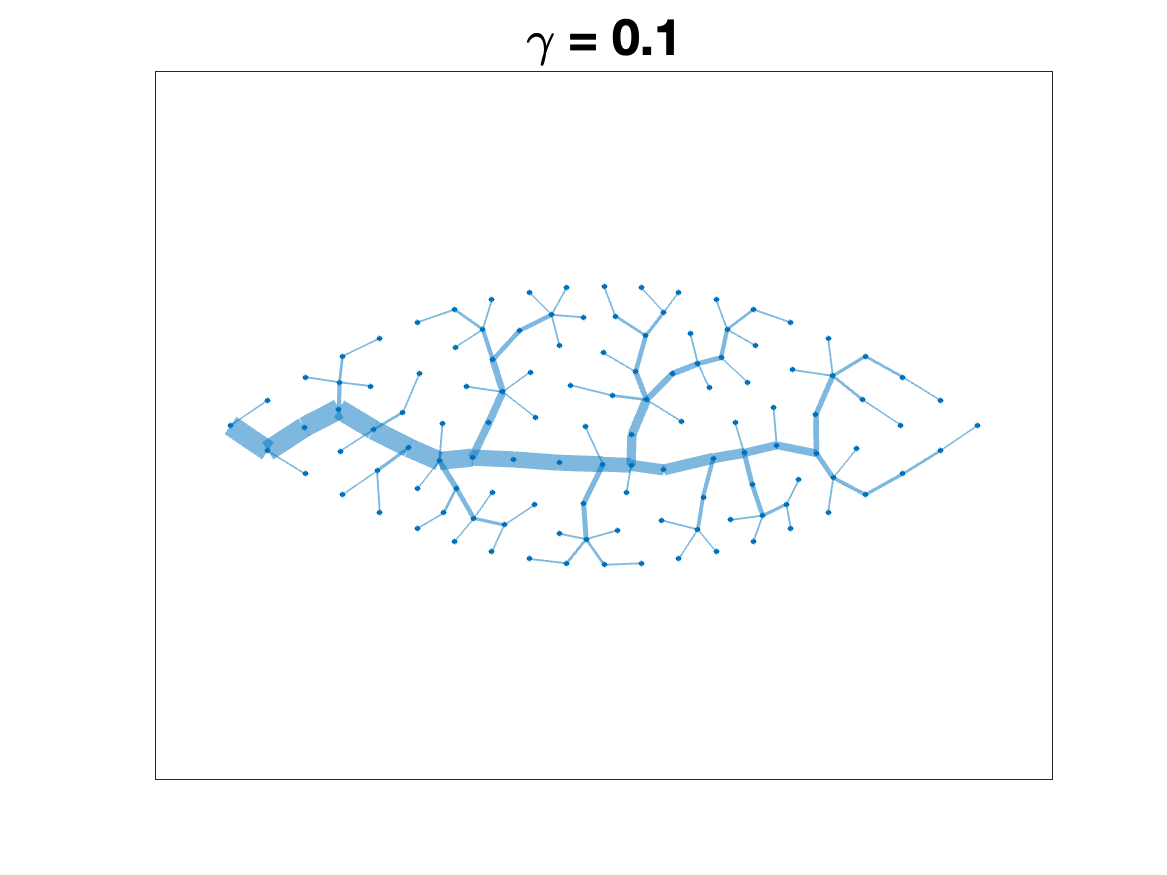}
 \includegraphics[width=.45\textwidth]{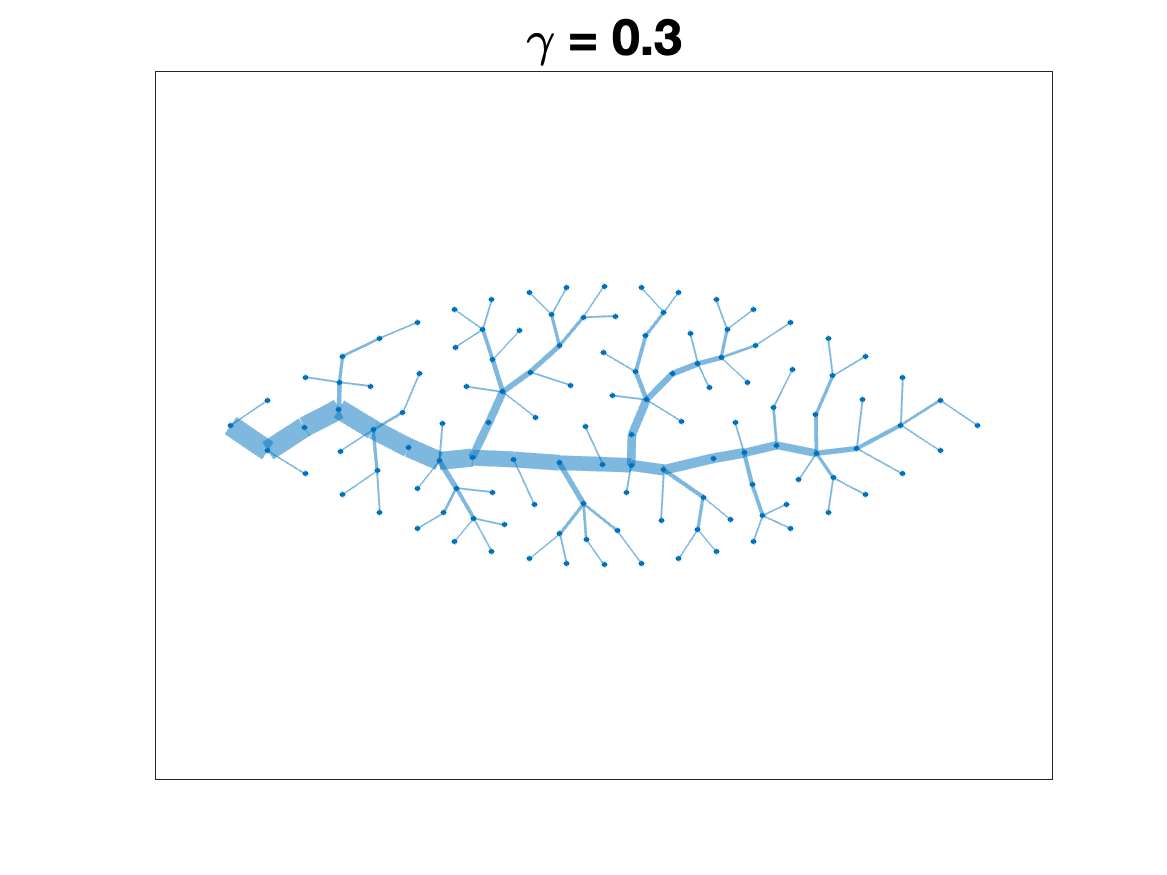} \\
 \includegraphics[width=.45\textwidth]{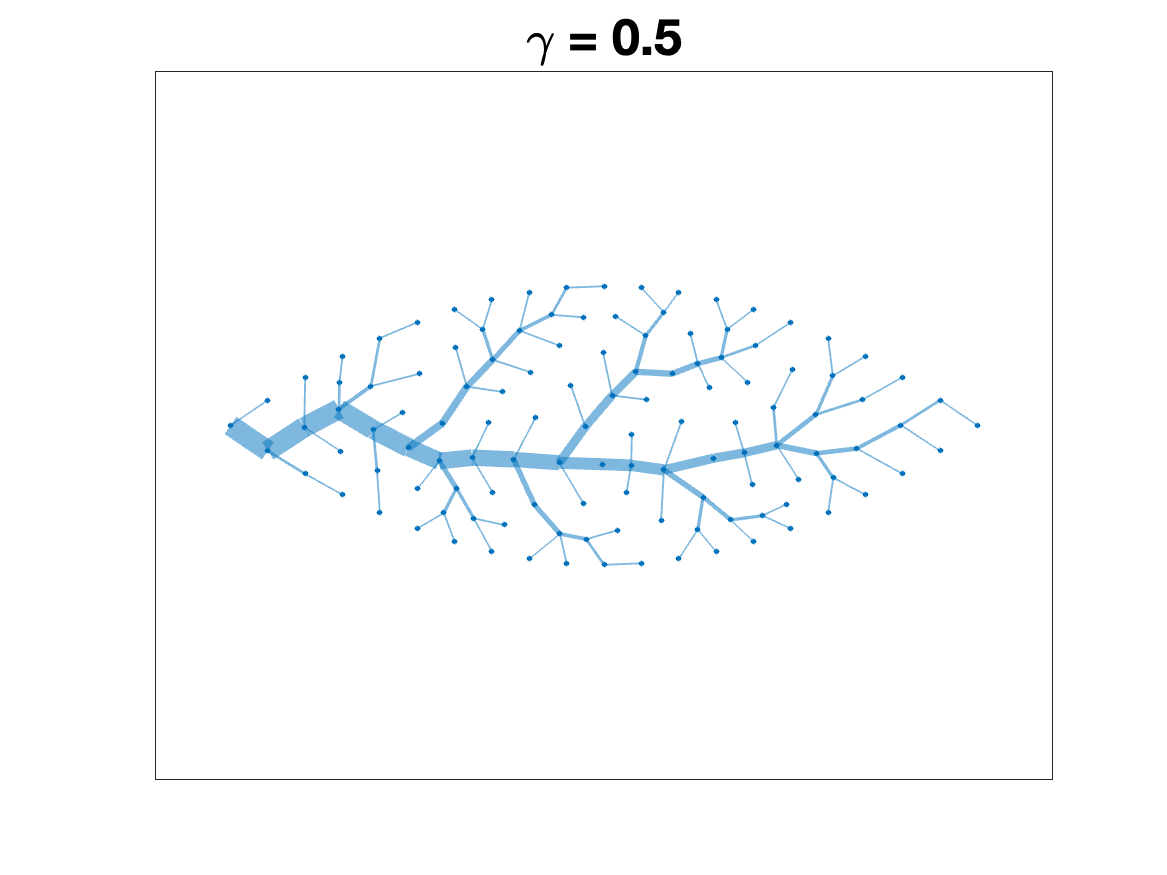} 
 \includegraphics[width=.45\textwidth]{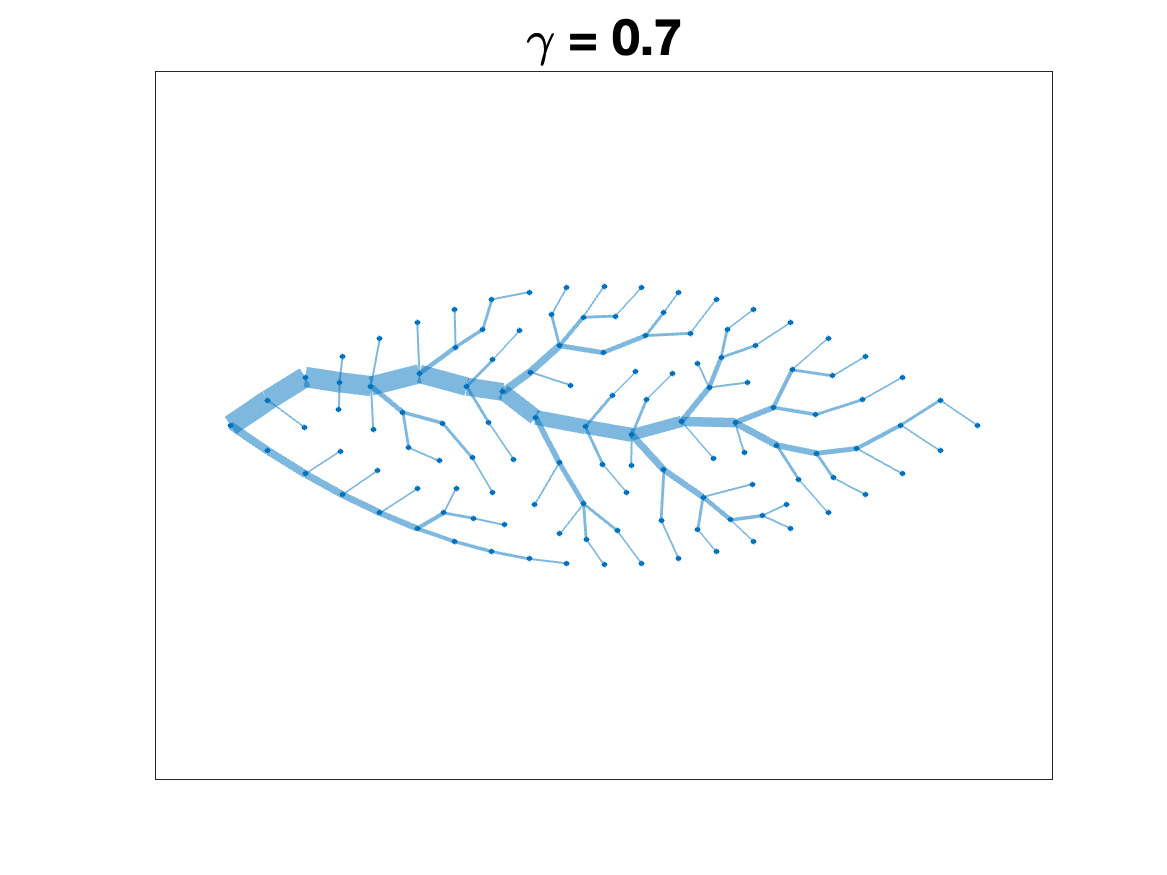} \\
 \includegraphics[width=.45\textwidth]{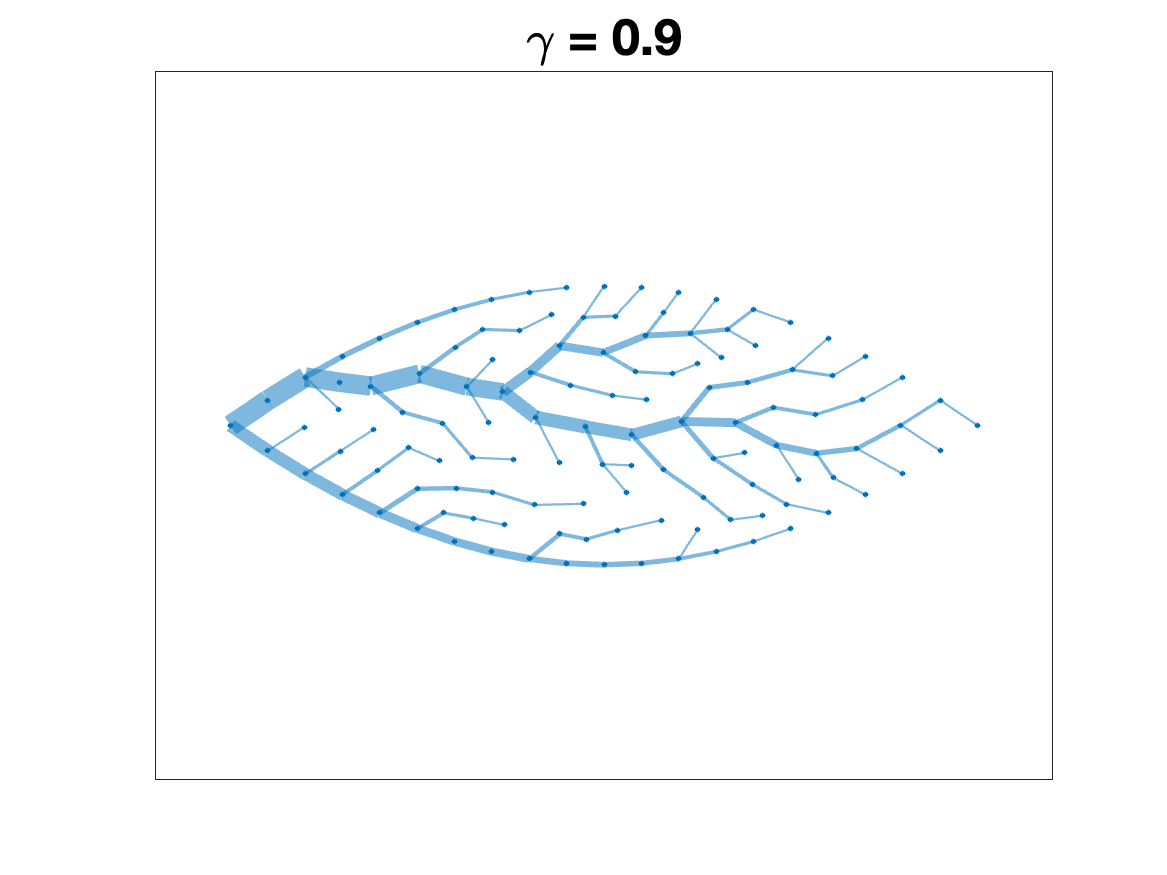}
 \includegraphics[width=.45\textwidth]{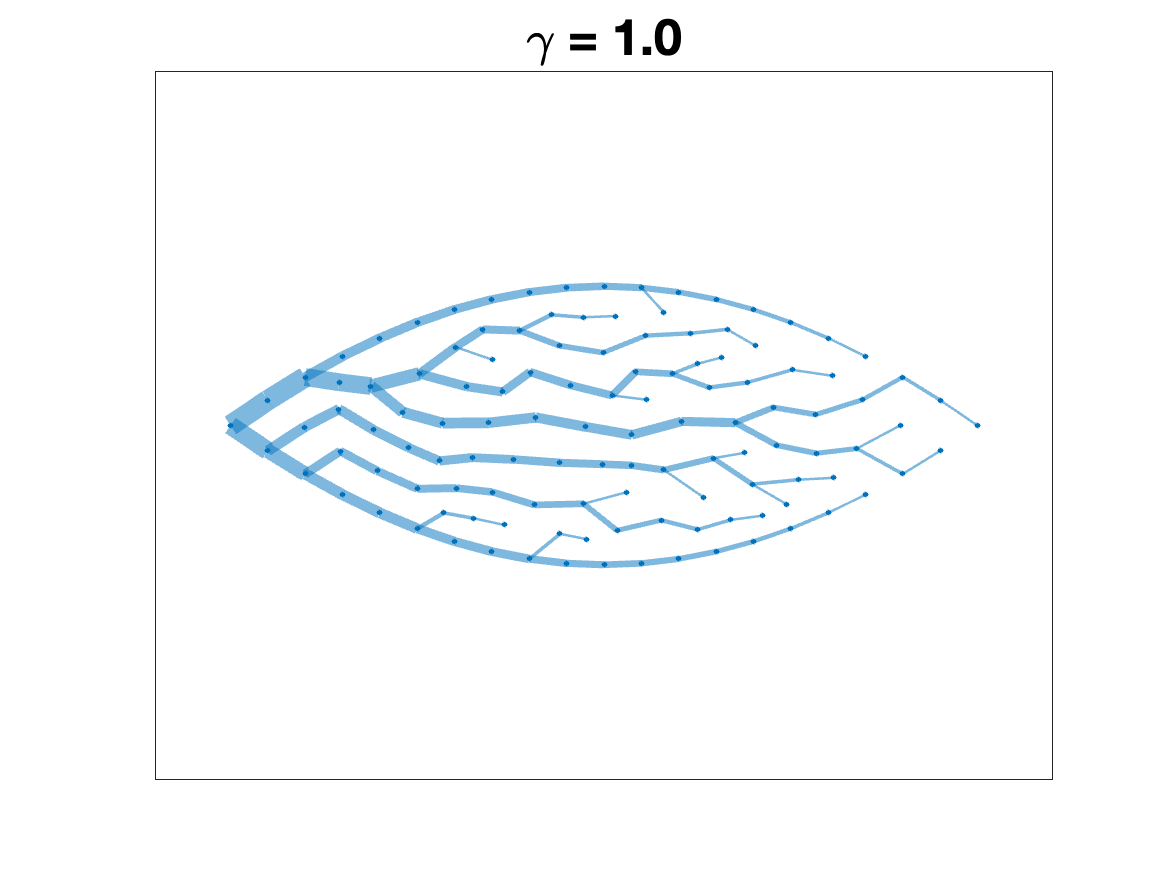} \\
 \caption{Results obtained by the discrete descent algorithm for the graph plotted in Fig. \ref{fig:leaf}
and  $\gamma\in\{0.1, 0.3, 0.5, 0.7, 0.9, 1.0\}$.
 The thickness of the line segments is proportional to the square root of the magnitude of the fluxes $|Q_{ij}|$ of the corresponding edge. Edges with $Q_{ij}=0$ are excluded from the plot.
 The direction of the flow is from the single source at the left-hand side to the right.
 \label{fig:res}}
\end{figure}

\begin{figure}[h!] \centering
 \includegraphics[width=.49\textwidth]{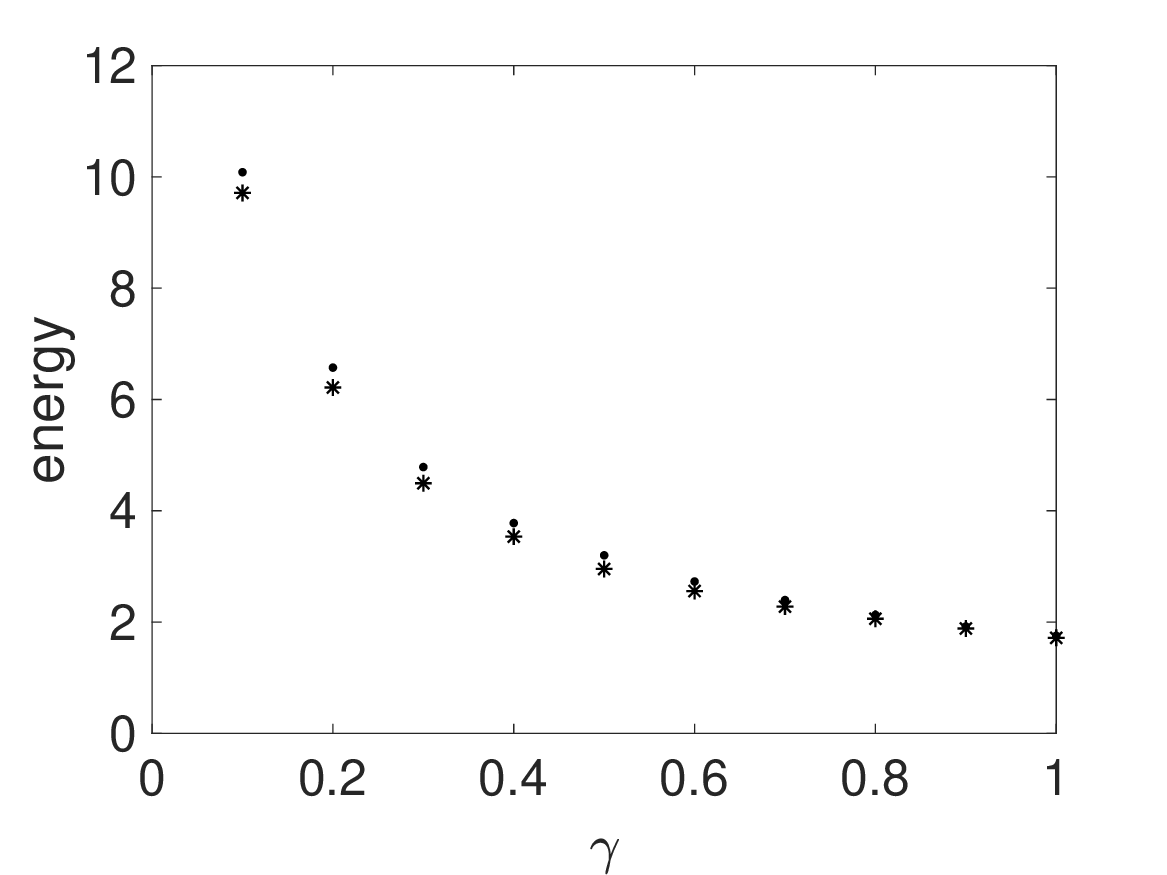}
 \includegraphics[width=.49\textwidth]{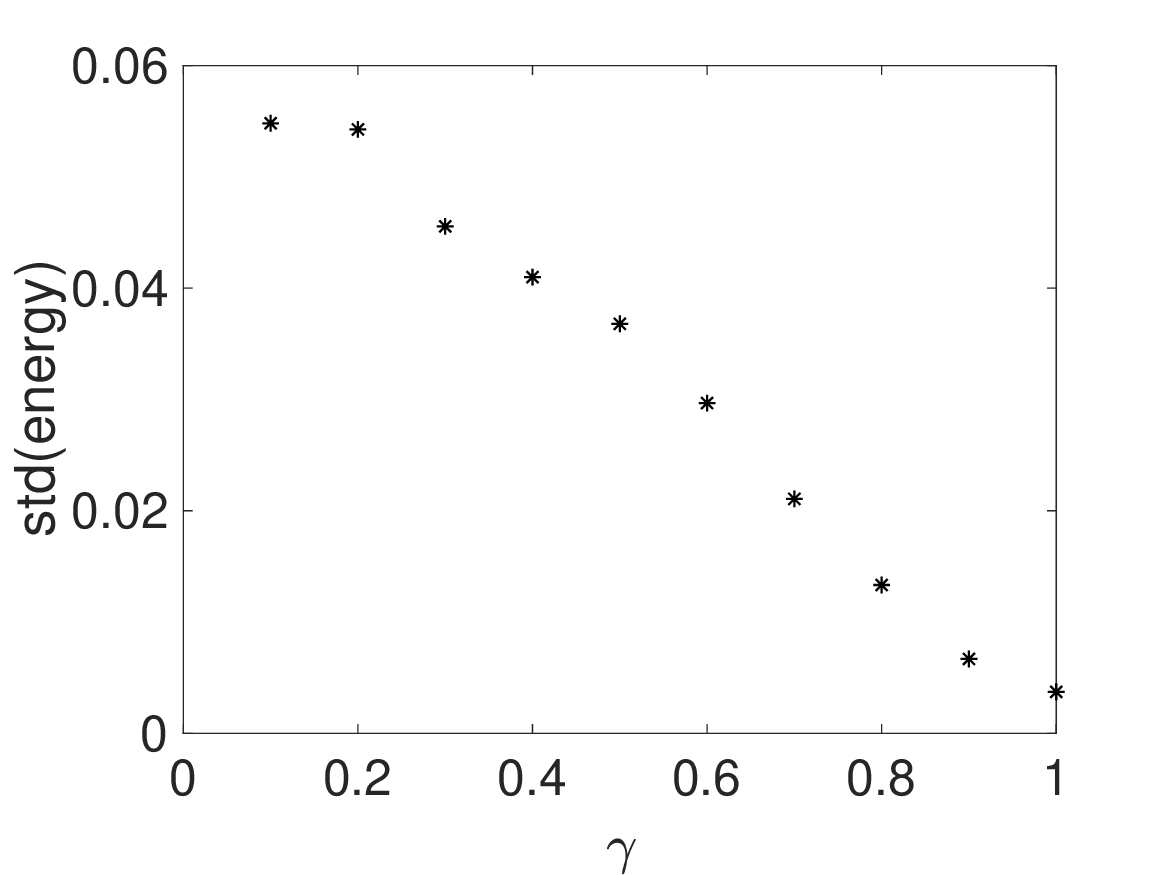} \\
 \caption{
 Minimal (star) and maximal (dot) values of the energy \eqref{eq:Enast}
 found in the $1000$ Monte-Carlo runs of the discrete descent algorithm (left panel)
 and their standard deviation (right panel).
 \label{fig:energy}}
\end{figure}

We carried out $1000$ Monte-Carlo runs of the discrete descent algorithm
for each of the values $\gamma\in\{0.1, 0.3, 0.5, 0.7, 0.9, 1.0\}$.
The graphs with the lowest value of the energy are presented in Fig. \ref{fig:res}.
We note that for $\gamma=1.0$ the output of the discrete descent algorithm
coincides with the global minimizer found by the gradient descent scheme, cf. \cite[Fig. 6]{HV23}.
To further illustrate the performance of the method, we plotted the maximal and minimal
values of the energy found during the $1000$ Monte-Carlo runs in left panel of Fig. \ref{fig:energy}
and the corresponding standard deviations in the right panel.
We observe that the performance of the method improves (i.e., the standard deviation decreases)
with increasing value of $\gamma$. This was to be expected since for larger values of $\gamma$
the energy functional becomes less non-convex (and it finally turns convex for $\gamma=1$).

Moreover, we investigated the degree of hierarchy of the (almost) optimal transportation structures
in dependence on the value of $\gamma\in(0,1]$.
As a measure of hierarchy we used the global reaching centrality (GRC) proposed in \cite{Vicsek}.
For unweighted directed networks it is defined as
\begin{equation}\label{eq:GRC}
   GRC = \frac{1}{|V|-1} \sum_{i\in\Vset} C_R^{\mathrm{max}} - C_R(i),
\end{equation}
where the local reaching centrality, $C_R(i)$, of node $i\in\Vset$ is the proportion of
all nodes in the graph that can be reached from node $i$ via outgoing edges,
and $C_R^{\mathrm{max}} := \max_{i\in\Vset} C_R(i)$.
In our case the direction of each edge is determined by the direction
of the corresponding flux, i.e., $i\to j$ if $Q_{ij}>0$ and vice versa, while edges with $Q_{ij}=0$
are, as always, treated as nonexistent.
In general, it would make sense to consider a weighted variant of the local reaching centrality,
defined as the sum of sources and sinks $S_j$ of nodes $j\in\Vset$ that can be reached from node $i$ via outgoing edges.
However, for the loop-free networks with only one source node and uniformly distributed sinks, as we consider here, this would lead to
the same result as the unweighted version.

\begin{figure}[h!] \centering
 \includegraphics[width=.49\textwidth]{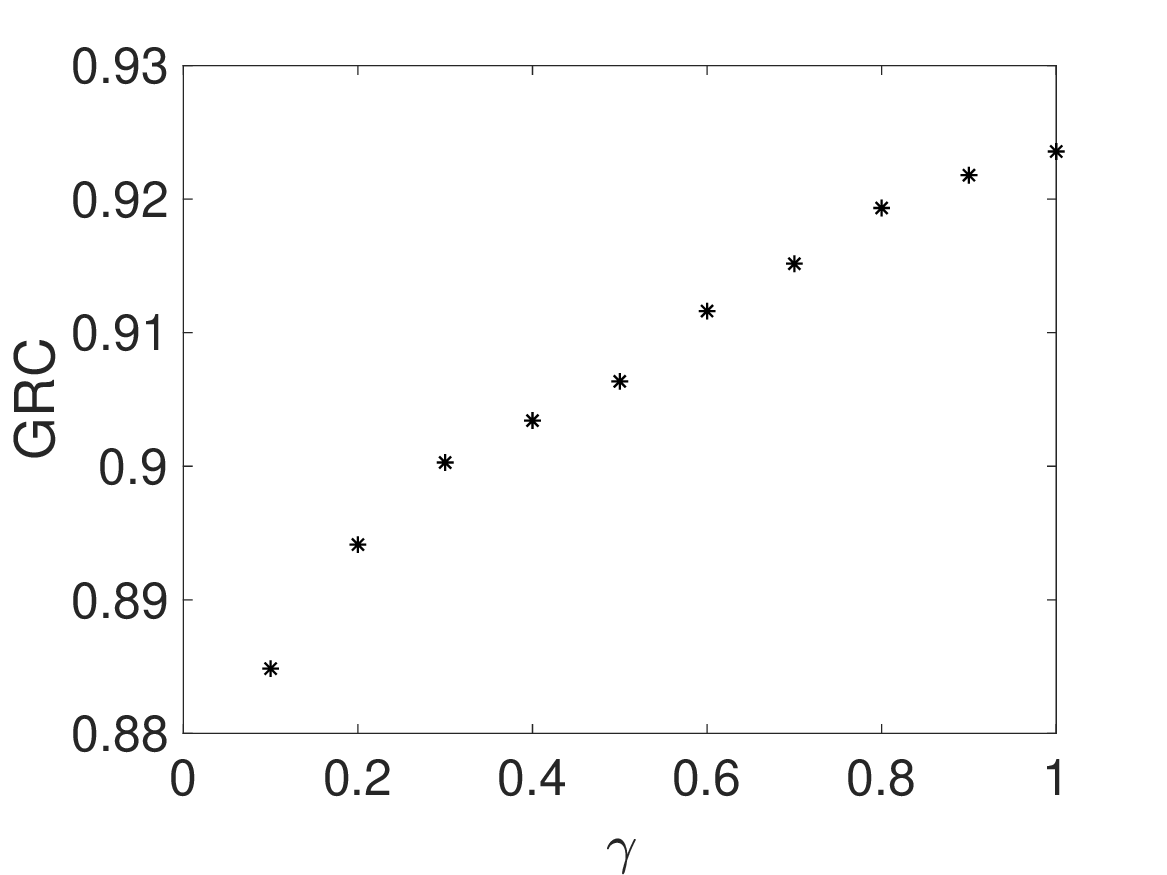}
 \includegraphics[width=.49\textwidth]{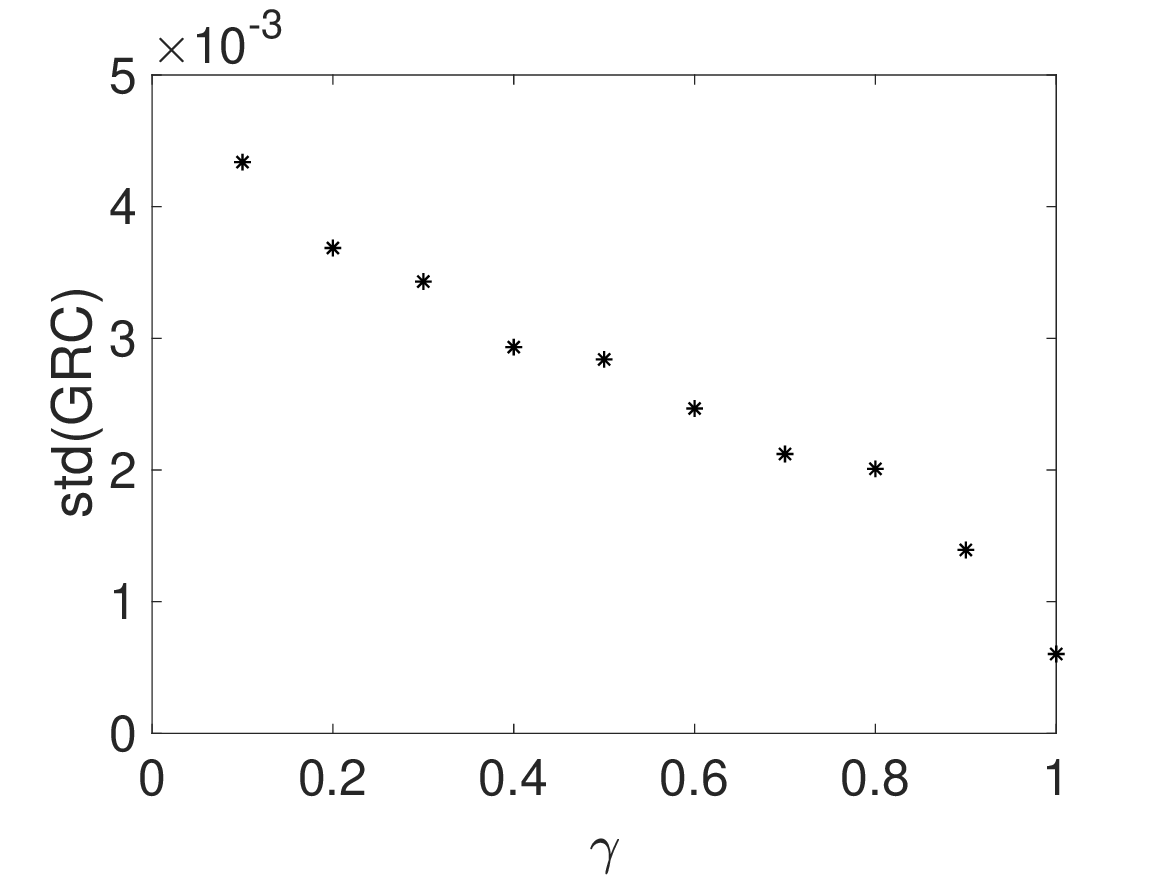} \\
 \caption{
 Values of the GRC coefficient \eqref{eq:GRC} for the networks with the lowest energies
found during the $1000$ runs of the Monte-Carlo simulation (left panel)
and their standard deviation (right panel).
 \label{fig:GRC}}
\end{figure}

The values of the GRC coefficient \eqref{eq:GRC} for the networks with the lowest energies
found during the $1000$ runs of the Monte-Carlo simulation are plotted in Fig. \ref{fig:GRC}, left panel.
We observe that the value of GRC monotonically increases with increasing $\gamma$,
which demonstrates that the value of the metabolic exponent controls the hierarchical organization
of the (almost) optimal transportation network.
To visualize the statistical properties of the GRC values produced by the Monte-Carlo simulation
we provide their standard deviation in the right panel of Fig. \ref{fig:GRC}.
We again observe that the standard deviation decreases with increasing value of the metabolic exponent $\gamma$.

\begin{thebibliography}{99}

\bibitem{bookchapter}
G. Albi, M. Burger, J. Haskovec, P. Markowich and M. Schlottbom, \emph{Continuum modelling of biological network formation},
in Active Particles Vol.I - Theory, Models, Applications, Series: Modelling and Simulation in Science and Technology,
N. Bellomo, P. Degond and E. Tamdor, eds., Birkh\"auser-Springer, Boston, 2017, pp. 1--48.
10.1007/978-3-319-49996-3

\bibitem{Ast1}
C. Astuto, D. Boffi, J. Haskovec, P. Markowich and G. Russo,
\emph{Asymmetry and condition number of an elliptic-parabolic system for biological network formation,}
Commun. Appl. Math. Comput. (2023).

\bibitem{Ast2}
C. Astuto, D. Boffi, J. Haskovec, P. Markowich and G. Russo,
\emph{Comparison of two aspects of a PDE model for biological network formation,}
Math. Comput. Appl. 27, 87 (2022).


\bibitem{Di}
D. Fang, S. Jin, P. Markowich and B. Perthame,
\emph{Implicit and Semi-implicit Numerical Schemes for the Gradient Flow of the Formation of Biological Transport Networks,}
SMAI Journal of Computational Mathematics 5 (2019), pp. 229--249.

\bibitem{Buchin-Schulz}
K. Buchin and A. Schulz, \emph{On the Number of Spanning Trees a Planar Graph Can Have}, in Algorithms -- ESA 2010. Lecture Notes in Computer Science 6346,
M. de Berg and U. Meyer, eds.,
Springer, Berlin, Heidelberg, 2010, pp. 110--121.

\bibitem{BHMR}
M. Burger, J. Haskovec, P. Markowich and H. Ranetbauer,
\emph{A mesoscopic model of biological transportation networks,}
Comm. Math. Sci., 17 (2019), pp. 1213--1234.


\bibitem{HKM}
J. Haskovec, L.-M. Kreusser and P. Markowich,
\emph{ODE and PDE based modeling of biological transportation networks,}
Comm. Math. Sci., 17 (2019), pp. 1235--1256.

\bibitem{HMP15} 
J. Haskovec, P. Markowich,  B. Perthame, \emph{Mathematical Analysis of a PDE System for Biological Network Formation,}
Communications in Partial Differential Equations, Vol. 40, No. 5, pp. 918--956, 2015.

\bibitem{HMPS}
J. Haskovec, P. Markowich, B. Perthame and M. Schlottbom,
\emph{Notes on a PDE system for biological network formation,}
Nonlinear analysis : theory, methods \& applications, 138 (2016), pp. 127--155.

\bibitem{HMP22}
J. Haskovec, P. Markowich and G. Pilli, \emph{Tensor PDE model of biological network formation,}
Comm. Math. Sci., 20 (2022), pp. 1173--1191.

\bibitem{HMP23}
J. Haskovec, P. Markowich and S. Portaro, \emph{Emergence of biological transportation networks as a self-regulated process,}
Discrete Contin. Dyn. Syst., 43 (2023), pp. 1499--1515.

\bibitem{HV23}
J. Haskovec and J. Vyb\'\i ral: \emph{Robust network formation with biological applications,}
preprint, arXiv:2311.17437 (2023).

\bibitem{Hu-Cai-13} D. Hu and D. Cai, \emph{Adaptation and optimization of biological transport networks,}
Phys. Rev. Lett., 111 (2013), pp. 138701.

\bibitem{Hu-Cai-19}
D. Hu and D. Cai, \emph{An optimization principle for initiation and adaptation of biological transport networks,}
Comm. Math. Sci. 17 (2019), pp. 1427--1436.

\bibitem{Laguna-Bohn-Jagla}
M. Laguna, S. Bohn and E. Jagla,
\emph{The Role of Elastic Stresses on Leaf Venation Morphogenesis,}
PLoS Comput Biol 4(4): e1000055 (2008).

\bibitem{Murray}
C. Murray, \emph{The physiological principle of minimum work. I. the vascular system and the cost of blood volume.}
Proc. Natl. Acad. Sci. USA, 12 (1926), pp. 207-214.

\bibitem{Runions}
A. Runions, M. Fuhrer, B. Lane, P. Federl, A.-G. Rolland-Lagan and P. Prusinkiewicz,
\emph{Modeling and visualization of leaf venation patterns,}
ACM Trans. Graph., 24 (2005), pp. 702--711.

\bibitem{Vicsek}
Mones E, Vicsek L, Vicsek T (2012), Hierarchy Measure for Complex Networks. PLoS ONE 7(3): e33799.

\end{thebibliography}
\end{document}